\documentclass[12pt]{iopart}

\usepackage{iopams}
\usepackage{graphicx}
\newtheorem{theorem}{Theorem}
\newtheorem{remark}{Remark}

\begin{document}

\title[Symmetries and exact solutions of the rotating shallow water equations]
{Symmetries and exact solutions of the rotating shallow water equations}

\author{Alexander Chesnokov}

\address{Lavrentyev Institute of Hydrodynamics SB RAS, Novosibirsk 630090, Russia}
\ead{chesnokov@hydro.nsc.ru}

\begin{abstract}
Lie symmetry analysis is applied to study the nonlinear rotating shallow water
equations. The 9-dimensional Lie algebra of point symmetries admitted by the
model is found. It is shown that the rotating shallow water equations are
related with the classical shallow water model with the change of variables.
The derived symmetries are used to generate new exact solutions of the rotating
shallow equations. In particular, a new class of time-periodic solutions with
quasi-closed particle trajectories is constructed and studied. The symmetry
reduction method is also used to obtain some invariant solutions of the model.
Examples of these solutions are presented with a brief physical interpretation.
\end{abstract}




\section{Introduction}
\eqnobysec

The article focuses on the Lie symmetry analysis of the rotating shallow water
(RSW) model. In Cartesian frame of reference the RSW equations have the
following form
\begin{equation}
\eqalign{
  u_t+uu_x+vu_y-fv+gh_x=0,  \cr
  v_t+uv_x+vv_y+fu+gh_y=0,  \cr
  h_t+(uh)_x+(vh)_y=0. } \label{eq:model}
\end{equation}
Here $(u,v)$ is the fluid velocity, $h$ is the free surface height over the
flat bottom, $f$ is the constant Coriolis parameter and $g$ is the constant
gravity acceleration. The RSW model (\ref{eq:model}) arises from the
three-dimensional rotating incompressible Euler equations (see \cite{Pedlosky},
pp. 59--63) under the following assumptions: a) the scale $H_0$ of the vertical
motion is much less than the scale $L_0$ of the horizontal motion, so that
$H_0/L_0<<1$; b) the fluid density is constant ($\rho={\rm const}$); c) the
external force is due to gravity, and pressure obeys the hydrostatic
approximation $p=pg(h-z)+p_0$ ($p_0$ is the constant pressure on the free
surface); d) the axis of rotation of the fluid coincides with the vertical
$z$-axis.

This nonlinear system of partial differential equations (PDEs) is widely-used
approximation for atmospheric and oceanic motions in the midlatitudes with
relatively large length and time scales. The model applies to phenomena that do
not depend substantially on temporal changes of the density stratification. In
spite of its simplicity, it contains all essential ingredients of atmosphere
and ocean dynamics at the synoptic scale. The main known mathematical results
and physical applications concerning of the RSW model (\ref{eq:model}) as well
as its linear theory, quasi-geostrophic limit and etc are presented in
\cite{Pedlosky}--\cite{Zeitlin}.

We recall here that the potential vorticity, defined by
$\Omega=(v_x-u_y+f)h^{-1}$, is conserved on particles:
$\Omega_t+u\Omega_x+v\Omega_y=0$. Moreover, there are two important for
applications non-dimensional parameters: Rossby number ${\rm Ro}=U_0/(fL_0)$
and Froude number ${\rm Fr}=U_0/\sqrt{gH_0}$ (here $U_0$ is the typical
velocity scale for the flow). In non-dimensional form of the model ${\rm
Ro}^{-1}$ and ${\rm Fr}^{-2}$ stand for $f$ and $g$, correspondingly. In this
study these parameters do not play a vital part. In fact, as result of obvious
scaling transformations the system (\ref{eq:model}) can be reduced to the
equivalent one with $f=1$ and $g=1$ (in case $f\neq 0$ and $g\neq 0$).

One of the methods for studying of PDEs is group analysis. This analytical
approach based on symmetries of differential equations was originally
introduced by Sophus Lie at the end-nineteenth century and further developed by
Ovsynnikov \cite{Ovs82}, Olver \cite{Olver} and others. For each system of PDEs
there is a symmetry group, that acts on the space of its independent and
dependent variables, leaving the form of the system unchanged. The classical
Lie symmetry analysis allows one to construct and study special types of
analytical solutions of nonlinear PDEs in terms of solutions of lower-dimension
equations. For modern description of the theory see, for example,
\cite{Ovs82}--\cite{Bluman}. Many applications of group analysis to PDEs are
collected in \cite{Ibr95}. The classical Lie method is an algorithmic procedure
for which many symbolic manipulation programs were designed \cite{Carminati,
Hereman}. This software became imperative in finding symmetries associated with
large systems of PDEs.

Exact solutions of nonlinear systems descriptive of fluid motions with moving
boundaries are uncommon even in the shallow water approximation. Freeman
\cite{Freeman} and Sachdev \cite{Sachdev80} have obtained solutions for gravity
wave motions under the assumption that the flow is a simple wave. It is
interesting to note that the similarity solution obtained in \cite{Freeman} may
be systematically generated by group methods \cite{Sachdev86}. In the latter
work, symmetry analysis was often used to obtain analytical solutions to the
nonlinear water-wave problems. Symmetries and exact solutions of a nonlinear
system that models the finite motion of a rotating shallow liquid contained in
a rigid basin were studied in \cite{Levi} (see also \cite{Rogers}). At present
there known many examples of exact (invariant and partial invariant) solutions
for nonlinear models of fluid dynamics \cite{Sachdev86}--\cite{Gol08} and
others. Amount of analytical solutions to the Euler equations govern the
compressible inviscid flow were obtained within the framework of Ovsynnikov's
``submodels programme'' \cite{Ovs94, Ovs99}.

In this paper we investigate Lie point symmetries and classes of exact
solutions of the RSW equations. In Section 2, we determine the 9-dimension Lie
algebra of infinitesimal transformations admitted by the RSW model. We also
prove that the derived Lie algebra of symmetries is isomorphic to the Lie
algebra of infinitesimal transformations admitted by (2+1)-dimensional shallow
water (SW) equations. This allows one to use its known optimal system of
subalgebras. Moreover, we find the change of variables which transforms the RSW
model to the SW equations. This is one of the principal results of the paper.
In Section 3, we derive and analyse finite transformations corresponding to the
nontrivial symmetries of the RSW model. On basis of these transformations, in
Section 4, we construct new time-periodic exact solutions of the RSW equations.
These solutions may be interpreted as pulsation of liquid volume under the
influence of gravity and Coriolis forces. In Section 5, using some
two-dimensional parameterized classes of the optimal system of subalgebras we
reduce the RSW model to the ordinary differential equations (ODEs) and
integrate them. In particular, we construct and study exact solutions
describing rotational-symmetric subcritical and supercritical regimes of the
stationary flow as well as various non-stationary regimes of spreading and
collapse of a liquid ring.

\section{Symmetries of the RSW equations and reduction to the SW model}

\subsection{Lie point symmetries of the RSW equations}

On the basis of the group analysis of the differential equations we shall find
the symmetry group of infinitesimal transformations associated with the RSW
equations (\ref{eq:model}). Following \cite{Ovs82}, let us define infinitesimal
generator $X$ and its first prolongation $Y$
\[ X=\xi^i({\bf x},{\bf u})\partial_{x^i}+\eta^j({\bf x},{\bf
   u})\partial_{u^j}, \quad
   Y=X+\zeta_i^j\partial_{u^j_i} \quad\quad (i,j=1,...,3). \]
Here we use the following notations for convenience
\[ {\bf x}=(x^1,x^2,x^3)=(t,x,y), \quad {\bf u}=(u^1,u^2,u^3)=(u,v,h) \]
\[ \zeta_i^j=D_i\eta^j-u^j_iD_i\xi^j, \quad
   u_i^j=\frac{\partial u^j}{\partial x^i}, \quad
   D_i=\frac{\partial}{\partial x^i}+u^j_i\frac{\partial}{\partial u^j}. \]
To determine the group of symmetries admitted by the system of equations
(\ref{eq:model}), we act on it with the first prolongation of generator $X$ and
go to the set of the system solutions. Consequently, we get a system of the
determining equations for the unknown functions $\xi^i({\bf x},{\bf u})$ and
$\eta^j({\bf x},{\bf u})$, that allows splitting in variables $u^j_i$. Let us
omit the bulky intermediate calculations and present the final result of the
symmetry group determination.

\begin{theorem} The symmetry group associated with the RSW equations
(\ref{eq:model}) is generated by the following vector fields:
\end{theorem}
\begin{equation}
 \eqalign{
   X_1=\partial_x, \quad X_2=\partial_y, \quad X_7=\partial_t, \cr
   X_3=\cos(ft)\partial_x-\sin(ft)\partial_y-f\sin(ft)\partial_u-
   f\cos(ft)\partial_v, \cr
   X_4=\sin(ft)\partial_x+\cos(ft)\partial_y+f\cos(ft)\partial_u-
   f\sin(ft)\partial_v, \cr
   X_5=-y\partial_x+x\partial_y-v\partial_u+u\partial_v, \cr
   X_6=x\partial_x+y\partial_y+u\partial_u+v\partial_v+ 2h\partial_h, }
 \label{eq:oper-f-1}
\end{equation}
\begin{eqnarray}
\fl  X_8= \cos(ft)\partial_t- \frac{f}{2}\Bigl(x\sin(ft)-
          y\cos(ft)\Bigr)\partial_x- \frac{f}{2}\Bigl(x\cos(ft)+
          y\sin(ft)\Bigr)\partial_y+ \nonumber\\
         +\frac{f}{2}\Bigl((u-f y)\sin(ft)+(v-f x)\cos(ft)\Bigr)\partial_u- \nonumber\\
         -\frac{f}{2}\Bigl((u+fy)\cos(ft)-(v+fx)\sin(ft)\Bigr)\partial_v+
         fh\sin(ft)\partial_h,
\label{eq:oper-f-2}
\end{eqnarray}
\begin{eqnarray}
\fl   X_9= \sin(ft)\partial_t+
        \frac{f}{2}\Bigl(x\cos(ft)+y\sin(ft)\Bigr)\partial_x-
        \frac{f}{2}\Bigl(x\sin(ft)-y\cos(ft)\Bigr)\partial_y- \cr
        -\frac{f}{2}\Bigl((u-fy)\cos(ft)-(v-fx)\sin(ft)\Bigr)\partial_u-\cr
        -\frac{f}{2}\Bigl((u+fy)\sin(ft)+(v+fx)\cos(ft)\Bigr)\partial_v-
        fh\cos(ft)\partial_h.
\label{eq:oper-f-3}
\end{eqnarray}

Thus, the RSW equations (\ref{eq:model}) are invariant under translations in
$t$, $x$ and $y$ ($X_7$, $X_1$ and $X_2$, correspondingly), helical rotations
with respect to $t$ ($X_4$ and $X_5$), rotation ($X_5$), scaling symmetry
($X_6$) and two complicated vector fields ($X_8$ and $X_9$) which will be
discussed below. The infinitesimal symmetries
(\ref{eq:oper-f-1})--(\ref{eq:oper-f-3}) form Lie algebra $L^f_9$.

Symmetry properties of the RSW equations in Lagrangian coordinates were studied
in \cite{Bila}. Transition to Lagrangian variables is a nonlocal change of
coordinate. Therefore Lie groups admitted by the model in Eulerian and
Lagrangian variables do not completely coincide. In particular, there are no
nontrivial symmetries defined by generators $X_8$ and $X_9$ in Lagrangian
variables.

As mentioned above, in \cite{Levi, Rogers} group analysis was applied to the
model of rotating shallow liquid in a rigid basin with the following geometry
\[ z=Z(x,y)=Ax^2+By^2, \quad A>0, \quad B>0. \]
It was shown that in the elliptic paraboloid case, $A\neq B$, the symmetry
algebra is spanned by six infinitesimal generators. In the circular paraboloid
case, corresponding to $A=B$, the symmetry algebra is larger --- namely,
9-dimensional. Note that these symmetries reduce to
(\ref{eq:oper-f-1})--(\ref{eq:oper-f-3}) for $A=B=0$.

\subsection{Properties of Lie algebra $L_9^f$}

Let $E$ be a given system of differential equations admitting a symmetry group
$G$. The basic property of solutions for $E$ is that any solution of the system
$E$ carried over by any transformation of the group $G$ to a certain solution
of the same system $E$. Therefore, two solutions of the system $E$ are said to
be essentially different with respect to $G$ if they can not transformed to
each other by any transformation of the group $G$. Thus, it is useful to
enumerate subgroups of the group $G$, which lead to essentially different
solutions. By virtue of the fact that there are certain correspondences between
the symmetry groups and Lie algebras, it is enough to construct an optimal
system of subalgebras $\Theta L$ (see \cite{Ovs82, Winter, Ovs93}), i.e. the
minimal set of subalgebras of $L$, which exhaust all the essentially different
invariant and partially invariant solutions to equations $E$.

In the present case there is no need to construct the optimal system of the
subalgebras for Lie algebra $L_9^f$ since we can use the results derived from
the symmetry analysis of the gas dynamic equations. In the case of zero
Coriolis parameter ($f=0$) the model under consideration reduces to the SW
equation (or two-dimensional polytropic gas dynamic equations with polytropic
exponent $\gamma=2$)
\begin{equation}
\eqalign{
  u_t+uu_x+vu_y+gh_x=0,  \cr
  v_t+uv_x+vv_y+gh_y=0,  \cr
  h_t+(uh)_x+(vh)_y=0. } \label{eq:SW}
\end{equation}
The Lie algebra $L_9$ of infinitesimal symmetries of the SW equations is
spanned by the following vector fields \cite{Ovs82, Ibr95}
\begin{equation*}
   \eqalign{
   Z_1=\partial_x, \ Z_2=\partial_y, \ Z_3=t\partial_x+\partial_u, \
   Z_4=t\partial_y+\partial_v, \cr
   Z_5=-y\partial_x+x\partial_y-v\partial_u+u\partial_v, \cr
   Z_6=x\partial_x+y\partial_y+u\partial_u+v\partial_v+2h\partial_h, \
   Z_7=\partial_t, \cr
   Z_8=t^2\partial_t+tx\partial_x+ty\partial_y+(x-tu)\partial_u+
   (y-tv)\partial_v-2th\partial_h, \cr
   Z_9=2t\partial_t+x\partial_x+y\partial_y-u\partial_u-v\partial_v-
   2h\partial_h. }
\end{equation*}
For this Lie algebra $L_9$, the optimal system of subalgebras $\Theta L_9$,
which includes 179 parameterized classes is presented in \cite{Pavl}. By virtue
of the following theorem, parameterized classes of this optimal system of
subalgebras can be applied for obtaining invariant or partially invariant
solutions of the RSW equations (\ref{eq:model}).

\begin{theorem}
Lie algebras $L_9(Z_1,...,Z_9)$ and $L_9^f(X_1,...,X_9)$ are isomorphic.
\end{theorem}

{\sf Proof}. Let us redefine the basis of generators of the Lie algebra
$L_9^f$ so that
\begin{equation}
  \eqalign{
   Y_1=X_2-X_4, \quad Y_2=X_3-X_1, \quad Y_3=X_1+X_3, \quad
   Y_4=X_2+X_4, \cr
   Y_5=X_5, \quad Y_6=X_6, \quad
   Y_7=\frac{1}{f}\Bigl(X_7-\frac{f}{2} X_5-X_8\Bigr), \cr
   Y_8=\frac{1}{f}\Bigl(X_7-\frac{f}{2} X_5+X_8\Bigr), \quad
   Y_9=-\frac{2}{f}X_9.
   }  \label{eq:oper-f-change}
\end{equation}
The commutation relations for the Lie algebra $L_9^f$ in this basis are
displayed in table \ref{tab1}. The basis is a canonical one chosen so that the
Levi decomposition is immediately apparent from the commutation table. Thus
$L_9^f=R\oplus N$, where $R=\{Y_1,...,Y_6\}$ and $N=\{Y_7,Y_8,Y_9\}$ constitute
the maximal solvable ideal (radical) and the simple Lie algebra $sl(2)$. The
nilradical of $L_9^f$ is abelian and is generated by $\{Y_1,Y_2,Y_3,Y_4\}$.

It is easy to verify that the commutation relations for the Lie algebra
$L_9(Z_1,...,Z_9)$ are defined by the same table of commutators (where symbols
$Y_k$ should be replaced by $Z_k$). The table of commutators for Lie algebras
$L_9(Z_1,...,Z_9)$ and $L^f_9(Y_1,...,Y_9)$ are coincide completely.
Consequently, Lie algebras $L_9$ and $L_9^f$ are isomorphic. $\Box$

\begin{table}
\begin{center}
\caption{\label{tab1} Table of commutators for Lie algebra
$L_9^f(Y_1,...,Y_9)$.}
\[\begin{array}{c|ccccccccc}
      & Y_1  & Y_2  & Y_3  & Y_4  & Y_5  & Y_6 & Y_7   & Y_8  & Y_9   \\ \hline
  Y_1 & 0    & 0    & 0    & 0    & Y_2  & Y_1 & 0     & Y_3  & Y_1   \\
  Y_2 & 0    & 0    & 0    & 0    & -Y_1 & Y_2 & 0     & Y_4  & Y_2   \\
  Y_3 & 0    & 0    & 0    & 0    & Y_4  & Y_3 & -Y_1  & 0    & -Y_3  \\
  Y_4 & 0    & 0    & 0    & 0    & -Y_3 & Y_4 & -Y_2  & 0    & -Y_4  \\
  Y_5 & -Y_2 & Y_1  & -Y_4 & Y_3  & 0    & 0   & 0     & 0    & 0     \\
  Y_6 & -Y_1 & -Y_2 & -Y_3 & -Y_4   & 0    & 0   & 0     & 0    & 0     \\
  Y_7 & 0    & 0    & Y_1  & Y_2  & 0    & 0   & 0     & Y_9  & 2Y_7  \\
  Y_8 & -Y_3 & -Y_4 & 0    & 0    & 0    & 0   & -Y_9  & 0    & -2Y_8 \\
  Y_9 & -Y_1 & -Y_2 & Y_3  & Y_4  & 0    & 0   & -2Y_7 & 2Y_8 & 0     \\
\end{array}\]
\end{center}
\end{table}

Thus, to construct invariant or partially invariant solutions to the RSW
equations (\ref{eq:model}) we can use generators $Y_k$ defined by formulas
(\ref{eq:oper-f-1})--(\ref{eq:oper-f-3}) and (\ref{eq:oper-f-change}) as well
as the parameterized classes of the optimal system of the  subalgebras
introduced in \cite{Pavl}.

\subsection{Reduction to the SW model}

In view of the previous theorem arises important question. Is there change of
variables such that each generator $Y_k$ $(k=1,...,9)$ reduces to generator
$Z_k$? In the case of existence one it is necessary to verify if the RSW
equations (\ref{eq:model}) reducible to the SW equations (\ref{eq:SW}) under
this transformation. The answer on the question is given by the following
theorem, which can be easy proved by straightforward calculations.

\begin{theorem}
The models (\ref{eq:model}) and (\ref{eq:SW}) as well as their symmetries are
related by the following transformation
\begin{equation}
 \fl
  \eqalign{
  t'=-\frac{1}{f}\cot\Bigl(\frac{ft}{2}\Bigr), \quad
  x'=-\frac{1}{2}\Bigl(x\cot\Bigl(\frac{ft}{2}\Bigr)-y\Bigr), \quad
  y'=-\frac{1}{2}\Bigl(x+y\cot\Bigl(\frac{ft}{2}\Bigr)\Bigr), \cr
  u'=-\frac{1}{2}\Bigl(u\sin(ft)-v(1-\cos(ft))-fx\Bigr), \cr
  v'=-\frac{1}{2}\Bigl(u(1-\cos(ft))+v\sin(ft)-fy\Bigr), \quad
  h'=\frac{h}{2}\Bigl(1-\cos(ft)\Bigr).
  } \label{eq:equiv-tr}
\end{equation}
\end{theorem}

In fact, let us suppose that a set of functions $u(t,x,y)$, $u(t,x,y)$,
$h(t,x,y)$ satisfies to the RSW equations (\ref{eq:model}). Then the immediate
calculations show that the following set of functions $u'(t',x',y')$,
$v'(t',x',y')$, $h'(t',x',y')$, given by formulae (\ref{eq:equiv-tr}),
satisfies to the SW equations (\ref{eq:SW}). Thus, the models (\ref{eq:model})
and (\ref{eq:SW}) are related by the change of variables (\ref{eq:equiv-tr}).

It is easy to clear up that each generator $Y_k$ in the variables with prime
(\ref{eq:equiv-tr}) takes form $Z_k$. More precisely, $Y'_k=Z_k$
($k=1,2,5,6,9$), $Y'_k=f Z_k$ ($k=3,4,8$), $Y'_7=f^{-1}Z_7$. Using scaling it
is possible to fix formulae (\ref{eq:equiv-tr}) (or reduce (\ref{eq:model}) to
the equivalent one with $f=1$) in order to $Y_k'=Z_k$ for all $k=1,...,9$
without coefficient $f$. However, this form is more convenient as it conserve
physical dimension of the variables.

Note that the mapping (\ref{eq:equiv-tr}) transforms any solution of the SW
equations (\ref{eq:SW}) defined for all time moments $(-\infty<t<\infty)$ to
the solution of the RSW equations (\ref{eq:model}) defined only for bounded
time interval ($0<t<2\pi/f$) or, may be, to the stationary solution. This
change of variables allows one to construct solutions of the RSW model using
known solutions of the SW equations and vice versa.

For example, let us consider constant solution of the SW equations
(\ref{eq:SW}): $u=u_0$, $v=v_0$, $h=h_0$. Obviously, in this case trajectories
of fluid particles are lines. According to (\ref{eq:equiv-tr}) we obtain the
following solution of the RSW model (\ref{eq:model})
\begin{equation}
  \eqalign{
    u=u_0\cot\Bigl(\frac{ft}{2}\Bigr)-v_0+
    \frac{f}{2}\Bigl(x\cot\Bigl(\frac{ft}{2}\Bigr)+y\Bigr), \cr
    v=u_0-v_0\cot\Bigl(\frac{ft}{2}\Bigr)-
    \frac{f}{2}\Bigl(x-y\cot\Bigl(\frac{ft}{2}\Bigr)\Bigr), \quad
    h=\frac{2h_0}{1-\cos(ft)}.
  }  \label{eq:const-SW}
\end{equation}
Solution (\ref{eq:const-SW}) is defined over the finite time interval
$0<t<2\pi/f$. Note that particle trajectories have the following form (after
time elimination)
\[ (x-A)^2+(y-B)^2=R^2, \]
where
\[ \fl A=x_0/2+u_0/f, \quad B=y_0/2+v_0/f, \quad
   R^2=(-x_0/2+u_0/f)^2+(y_0/2-v_0/f)^2. \]
Here $x_0$ and $y_0$ are coordinates of the particle at $t=\pi/f$. Thus, each
fluid particle moves around a circumference. However, these trajectories do not
close because solution becomes invalid as time $t$ tends to $0$ or $2\pi/f$.
This solution describes essentially two-dimension fluid flow. In the case
$u_0^2+v_0^2=0$ it reduces to the rotationally symmetric solution in the polar
frame of reference (\ref{eq:polar})
\[ U=\frac{fr}{2}\cot\Bigl(\frac{ft}{2}\Bigr), \quad V=-\frac{fr}{2}, \quad
   h=\frac{2h_0}{1-\cos(ft)}. \]

Inversely, let us take solution $u=0$, $v=0$, $h=h_0$ of the RSW equations
(\ref{eq:model}), corresponding to the state of rest. Transformation
(\ref{eq:equiv-tr}) yields the following rotationally-symmetric solution of the
SW equations (\ref{eq:SW})
\[ U=\frac{f^2tr}{1+f^2t^2}, \quad V=\frac{fr}{1+f^2t^2}, \quad
   h=\frac{h_0}{1+f^2t^2} \]
(primes are omitted here).

Using analogy with gas dynamic, we can conclude that these solutions belong to
the class of ``barochronous'' \cite{Chup97, Chup99} motions (as depth $h$
depends only on time $t$).

\section{Finite transformations corresponding to the ``complicate'' generators}

It is well known that generator $X=\xi^i\partial_{x^i}$ can be associated with
one-parameter group of transformation $G_1$ given by the finite relations of
the form $\bar{x}^i=\bar{x}^i(x^1,...,x^N,a)$, where $a$ is a real-valued
parameter. To determine these relations we shall solve the system of Lie
equations
\begin{equation}
  \frac{\partial\bar{x}^i}{\partial a}=\xi^i(\bar{x}^1,...,\bar{x}^N), \quad
  \bar{x}^i|_{a=0}=x^i \quad (i=1,...,N).
 \label{eq:Lie}
\end{equation}

Let us find the finite transformations corresponding to the nontrivial
infinitesimal generators $Y_7$, $Y_8$ and $Y_9$. For the sake of convenience we
shall use polar coordinates $(r,\theta)$, radial $U$ and circular $V$ velocity
vector components related to $(x,y)$ coordinates and $(u,v)$ velocities, so
that
\begin{equation}
 \fl
  x=r\cos\theta, \quad y=r\sin\theta, \quad
  u=U\cos\theta-V\sin\theta, \quad v=U\sin\theta+V\cos\theta.
 \label{eq:polar}
\end{equation}
In these variables admissible symmetries $Y_7$, $Y_8$, $Y_9$ take the following
form
\begin{eqnarray*}
 \fl  \hat{Y}_7= f^{-1}(1-\cos(ft))\partial_t+2^{-1}r\sin(ft)\partial_r-2^{-1}
      (1-\cos(ft))\partial_\theta-\cr
      -2^{-1}(U\sin(ft)-fr\cos(ft))\partial_U-
      2^{-1}(V+fr)\sin(ft)\partial_V-h\sin(ft)\partial_h,
\end{eqnarray*}
\begin{eqnarray*}
 \fl  \hat{Y}_8=-\hat{Y}_7+2f^{-1}\partial_t-\partial_\theta,
\end{eqnarray*}
\begin{eqnarray*}
 \fl  \hat{Y}_9= -2f^{-1}\sin(ft)\partial_t-r\cos(ft)\partial_r+
      \sin(ft)\partial_\theta+\cr
      +(U\cos(ft)+fr\sin(ft))\partial_U+
      (V+fr)\cos(ft)\partial_V+2h\cos(ft)\partial_h.
\end{eqnarray*}
Integration of the Lie equations (\ref{eq:Lie}) associated with generator
$\hat{Y}_9$ gives the following result ($t\neq (2n+1)\pi/f$, $n$ --- integer):
\begin{equation}
 \eqalign{
   \bar{t}=\frac{2}{f}\arctan(\alpha\tau)+\chi(t), \quad
   \bar{r}=r\sqrt{\frac{\alpha (1+\tau^2)}{1+\alpha^2\tau^2}}, \cr
   \bar{\theta}=\theta+\arctan(\tau)-\arctan(\alpha\tau), \cr
   \bar{U}=\left( U-\frac{fr}{2}\frac{(\alpha^2-1)\tau}{1+\alpha^2\tau^2}\right)
      \sqrt{\frac{1+\alpha^2\tau^2}{\alpha (1+\tau^2)}}, \cr
   \bar{V}=\left( V+\frac{fr}{2}\frac{(\alpha-1)(\alpha\tau^2-1)}{1+\alpha^2\tau^2}\right)
      \sqrt{\frac{1+\alpha^2\tau^2}{\alpha (1+\tau^2)}}, \quad
   \bar{h}=h\frac{1+\alpha^2\tau^2}{\alpha (1+\tau^2)}. } \label{eq:Lie-sol-1}
\end{equation}
Here $\tau=\tan(\frac{ft}{2})$, \ $\alpha=\exp(-2a)$, \ $\chi(t)=2\pi k/f$,
where integer $k$ is such that value $t$ specified by initial condition at
$a=0$ belongs to interval $((2k-1)\pi/f,(2k+1)\pi/f)$. In the case
$t=(2n+1)\pi/f$ solution of the equations (\ref{eq:Lie}) have the following
form
\begin{equation}
 \eqalign{
   \bar{t}=t=(2n+1)\pi/f, \quad \bar{r}=r/\sqrt{\alpha}, \quad
   \bar{\theta}=\theta, \cr
   \bar{U}=U\sqrt{\alpha}, \quad
   \bar{V}=\left(V+(\alpha-1)(2\alpha)^{-1}fr\right)\sqrt{\alpha}, \quad
   \bar{h}=\alpha h.} \label{eq:Lie-sol-2}
 \end{equation}

Let us show that functions (\ref{eq:Lie-sol-1}) can be defined at points
$t=t_*=(2n+1)\pi/f$ by continuity in compliance with (\ref{eq:Lie-sol-2}). The
derived mapping will be continuously differentiable with respect to all its
arguments. Let us investigate one-sided limits in formula (\ref{eq:Lie-sol-1})
at $t\to t_*\mp 0$ (the remaining variables are fixed). In this case
$\tau(t)=\tan(\frac{ft}{2})\to\pm\infty$. Discontinuous at point $t=t_*$
function $\frac{2}{f}\arctan(\alpha\tau(t))$ has the following limits: $\pi/f$
on the left-hand side and $-\pi/f$ on the right-hand side of the discontinuity.
Function $\chi(t)$ is piecewise constant and possesses the value of $2\pi n/f$
on the left-hand side and $2\pi (n+1)/f$ on the right-hand side of point
$t=t_*$. Thus, function $\bar{t}(t,\alpha)= \frac{2}{f} \arctan(\alpha\tau)+
\chi(t)$ with $t\to t_*\mp 0$ has one-sided limits which coincide and are equal
to $(2n+1)\pi/f$. This agrees with the first formula in (\ref{eq:Lie-sol-2}).
Correspondence of the remain functions in formulae (\ref{eq:Lie-sol-1}) and
(\ref{eq:Lie-sol-2}) is easily established by calculation of their one-side
limits at $t=t_*$. Moreover, straightforward analysis shows that the mapping
(\ref{eq:Lie-sol-1}) redefined at points $t=t_*=(2n+1)\pi/f$ in compliance with
formulae (\ref{eq:Lie-sol-2}) is continuously differentiable with respect to
all the arguments.

Similarly we calculate finite transformations corresponding to the generator
$\hat{Y}_8$:
\begin{equation}
 \eqalign{
   \bar{t}=\frac{2}{f}\arctan(\tau+a)+\chi(t), \quad
   \bar{r}=r\sqrt{\frac{\tau^2+1}{(\tau+a)^2+1}}, \cr
           \bar{\theta}=\theta+\arctan(\tau)-\arctan(\tau+a), \cr
   \bar{U}=\left( U+\frac{fr}{2} \frac{(\tau^2 +a\tau-1)a}{(\tau+a)^2+1} \right)
           \sqrt{\frac{(\tau+a)^2+1}{\tau^2+1}}, \cr
   \bar{V}=\left( V+\frac{fr}{2}\frac{(2\tau+a)a}{(\tau+a)^2+1}\right)
           \sqrt{\frac{(\tau+a)^2+1}{\tau^2+1}}, \quad
   \bar{h}=h\frac{(\tau+a)^2+1}{\tau^2+1}. } \label{eq:Lie-sol-3}
\end{equation}
Here $\tau$ and $\chi$ are the same functions as before.

Finite transformations corresponding to the generator $\hat{Y}_7$ have the same
form (\ref{eq:Lie-sol-3}), where $\sigma=-\cot(\frac{ft}{2})$ and $\chi_1(t)$
stand for $\tau$ and $\chi(t)$. Function $\chi_1(t)$ is equal to $(2k+1)\pi/f$,
where integer $k$ is such that value $t$ specified by initial condition in the
system (\ref{eq:Lie}) at $a=0$ belongs to interval $(2k\pi/f,2(k+1)\pi/f)$.

According to the general theory \cite{Ovs82, Olver} governing equations
(\ref{eq:model}) presented in the polar coordinates (\ref{eq:polar})
\begin{equation}
 \eqalign{
  \frac{\partial U}{\partial t}+ U\frac{\partial U}{\partial r}+
   \frac{V}{r} \frac{\partial U}{\partial \theta}-\frac{V^2}{r}-fV+
   g\frac{\partial h}{\partial r}=0, \cr
   \frac{\partial V}{\partial t}+ U\frac{\partial V}{\partial r}+
   \frac{V}{r} \frac{\partial V}{\partial \theta}+ \frac{UV}{r}+
   fU+\frac{g}{r}\frac{\partial h}{\partial \theta}=0, \cr
   \frac{\partial h}{\partial t}+
   \frac{1}{r}\frac{\partial(rUh)}{\partial r}+
   \frac{1}{r}\frac{\partial (Vh)}{\partial \theta}=0} \label{eq:model-pol}
\end{equation}
remain unchanged under the substitutions (\ref{eq:Lie-sol-1}) and
(\ref{eq:Lie-sol-3}) corresponding to the finite transformations of the
``complicate'' generators. This fact allows one to formulate the following
theorem based on the finite transformations (\ref{eq:Lie-sol-1}). Similar
theorem can be obtained using formulae (\ref{eq:Lie-sol-3}).

\begin{theorem}
If a set of functions
\begin{equation*}
  U=\bar{U}(t,r,\theta), \quad V=\bar{V}(t,r,\theta), \quad
  h=\bar{h}(t,r,\theta)
\end{equation*}
satisfies the system of equations (\ref{eq:model-pol}) then the following set
of functions
\begin{equation}
 \eqalign{
   U(t,r,\theta)=\left( \bar{U}(\bar{t},\bar{r},\bar{\theta})+
     \frac{f\bar{r}}{2}\frac{(\alpha^2-1)\bar{\tau}}{\alpha^2+\bar{\tau}^2}\right)
     \sqrt{\frac{\alpha^2+\bar{\tau}^2}{\alpha (1+\bar{\tau}^2)}}, \cr
   V(t,r,\theta)=\left( \bar{V}(\bar{t},\bar{r},\bar{\theta})-
     \frac{f\bar{r}}{2}\frac{(\alpha-1)(\bar{\tau}^2-\alpha)}{\alpha^2+\bar{\tau}^2}\right)
     \sqrt{\frac{\alpha^2+\bar{\tau}^2}{\alpha (1+\bar{\tau}^2)}}, \cr
   h(t,r,\theta)=\frac{\alpha^2+\bar{\tau}^2}{\alpha (1+\bar{\tau}^2)}
     \bar{h}(\bar{t},\bar{r},\bar{\theta})}   \label{eq:new-sol}
\end{equation}
also satisfies the same system. Here
\begin{equation}
 \eqalign{
  \tau=\tan\left(\frac{ft}{2}\right), \quad
  \bar{\tau}=\tan\left(\frac{f\bar{t}}{2}\right), \quad
  \bar{t}=\frac{2}{f}\arctan(\alpha\tau)+\chi(t), \cr
  \bar{r}=r\sqrt{\frac{\alpha (1+\tau^2)}{1+\alpha^2\tau^2}}, \quad
  \bar{\theta}=\theta+\arctan(\tau)-\arctan(\alpha\tau), }  \label{eq:new-var}
\end{equation}
where $\alpha$ is an arbitrary positive constant, $\chi(t)$ is the piecewise
constant function which is equal to $2\pi k/f$ on the interval $t \in
((2k-1)\pi/f,(2k+1)\pi/f)$, $k$ --- integer. We assume that $\bar{t}=t$, \
$\bar{r}=r/\sqrt{\alpha}$, \ $\bar{\theta}=\theta$ at points $t=(2k+1)\pi/f$.
\end{theorem}

\section{Time-periodic exact solutions of the RSW equations}

The above-formulated theorem makes it possible to obtain new exact solutions of
the RSW model by using known ones. Similar approach was applied to the gas
dynamic equations in \cite{Nik, Oliveri}. Obviously, equations
(\ref{eq:model-pol}) have the following class of stationary
rotationally-symmetric solutions
\begin{equation}
   U=\bar{U}=0, \quad V=\bar{V}(r), \quad
   h=\bar{h}(r)=\frac{1}{g}\int\limits_0^r
   \left(\frac{\bar{V}^2}{r}+f\bar{V}\right)\,dr+h_0,
 \label{eq:st-sol-1}
\end{equation}
where $\bar{V}(r)$ is an arbitrary smooth function, $h_0$ is a positive
constant. According to the preceding theorem and formulae (\ref{eq:new-sol})
the functions
\begin{equation}
  \eqalign{
   U=\frac{fr}{2}\,\frac{(\alpha^2-1)\tau}{1+\alpha^2\tau^2}, \cr
   V=\sqrt{\frac{\alpha(1+\tau^2)}{1+\alpha^2\tau^2}}\,\bar{V}
   \left(r\sqrt{\frac{\alpha(1+\tau^2)}{1+\alpha^2\tau^2}}\right)-
   \frac{fr}{2}\, \frac{(\alpha-1)(\alpha\tau^2-1)}{1+\alpha^2\tau^2}, \cr
   h=\frac{\alpha(1+\tau^2)}{1+\alpha^2\tau^2}\,
   \bar{h}\left(r\sqrt{\frac{\alpha(1+\tau^2)}{1+\alpha^2\tau^2}}\right); \quad
   \tau=\tan\left(\frac{ft}{2}\right) }  \label{eq:periodic-sol}
\end{equation}
is the solution of the equations (\ref{eq:model-pol}). The principal
prerequisite to the existence of periodic solutions of the form
(\ref{eq:periodic-sol}) is the presence of nonzero Coriolis parameter $f$. Let
us consider in detail two typical solutions from class (\ref{eq:periodic-sol}).

\subsection{Example \#1: Pulsation of the liquid cylinder}

Let us consider solution from class (\ref{eq:st-sol-1}) corresponding to the
state of rest: $\bar{U}=0$, $\bar{V}=0$, and $\bar{h}=h_0$. By substitution of
these functions into (\ref{eq:periodic-sol}), the time-periodic exact solution
of the RSW equations (\ref{eq:model-pol})
\begin{equation}
 \eqalign{
  U=\frac{fr}{2}\,\frac{(\alpha^2-1)\tau}{1+\alpha^2\tau^2}, \quad
  V= -\frac{fr}{2}\, \frac{(\alpha-1)(\alpha\tau^2-1)}{1+\alpha^2\tau^2}, \cr
  h=\frac{\alpha(1+\tau^2)h_0}{1+\alpha^2\tau^2} }
 \label{eq:period-sol-1}
\end{equation}
is obtained. Here, as before, $\tau=\tan(\frac{ft}{2})$. Particle trajectories
are defined by the system of equations
\begin{equation}
   \frac{dr}{dt}=U, \quad \frac{d\theta}{dt}=\frac{V}{r}, \quad r|_{t=0}=r_0,
   \quad \theta|_{t=0}=\theta_0
  \label{eq:tr-eq}
\end{equation}
and in this case have the following form
\begin{equation}
   r(t)=r_0\sqrt{\frac{1+\alpha^2\tau^2(t)}{1+\tau^2(t)}}, \quad
   \theta(t)=\theta_0+\arctan(\alpha\tau(t))-\arctan(\tau(t)).
  \label{eq:tr-1}
\end{equation}
Formulae (\ref{eq:tr-1}) can be rewritten in a simpler form
\[ (x-A)^2+(y-B)^2=R^2, \]
where
\[ A=2^{-1}(\alpha+1)r_0\cos\theta_0, \quad
   B=2^{-1}(\alpha+1)r_0\sin\theta_0, \quad R=2^{-1}(\alpha-1)r_0. \]
Thus, each fluid particle moves round a circumference of radius $R$, which is
proportional to the initial distance of the particle from the origin of the
system of coordinates. All particles return to their initial positions at time
$T=2\pi/f$, that is period of solution (\ref{eq:period-sol-1}). The solution
can be interpreted as pulsation of the liquid cylinder with impermeable
boundary moving in accordance with the first formula (\ref{eq:tr-1}), where
$r_0$ is the initial cylinder boundary radius. The height of the liquid column,
depending only on time, is given by last formula in (\ref{eq:period-sol-1})
($\alpha h_0$ is fluid surface height at initial time $t=0$). We also note that
time-periodic solution (\ref{eq:period-sol-1}) is the solution with constant
potential vorticity $\Omega=f/h_0$.

\subsection{Example \#2: Pulsation of the liquid ``drop"}

Let us consider the following solution from class (\ref{eq:st-sol-1})
\begin{equation}
  \bar{U}=0, \quad \bar{V}(r)=l r^2, \quad \bar{h}(r)=\frac{l^2}{4g}
  \left(r^4+\frac{4f}{3l}r^3+\frac{f^4}{3l^4}\right),
  \label{eq:drop}
\end{equation}
where $l=-f^2\sqrt{\alpha/(12g)}<0$. In this case the fluid is localized on the
rotating plane within the cylinder $0\leq r\leq -f/l$ (we suppose, that $f>0$).
The depth of fluid $h$ is maximal at $r=0$ and vanishes at $r=-f/l$. According
to formulae (\ref{eq:periodic-sol}) and (\ref{eq:drop}) we arrive to the new
time-periodic exact solution of the RSW equations (\ref{eq:model-pol}), which
describes pulsation of liquid volume under the influence of gravitation and
rotation.

Integrating the equations of trajectories (\ref{eq:tr-eq}) with the functions
$U$, $V$ defined by (\ref{eq:periodic-sol}) and (\ref{eq:drop}), we obtain the
following result
\begin{equation}
 \eqalign{
   r(t)=r_0\sqrt{\frac{1+\alpha^2\tau^2(t)}{1+\tau^2(t)}},\cr
   \theta(t)=(2Cf^{-1}+1)\arctan(\alpha\tau(t))-
   \arctan(\tau(t))+ C\,\chi(t) +\theta_0. }
 \label{eq:tr-periodic}
\end{equation}
Here $C=lr_0\sqrt{\alpha}$; $\chi(t)$ is a piecewise constant function which
takes the value $2\pi k/f$ for $t$ from the interval
$((2k-1)\pi/f,(2k+1)\pi/f)$, where $k$ is an integer. Constants $r_0$ and
$\theta_0$ specify the position of the particle on the plane at initial time
$t=0$.

\begin{figure}[tbp]
\begin{center}
\resizebox{0.9\textwidth}{!}{\includegraphics{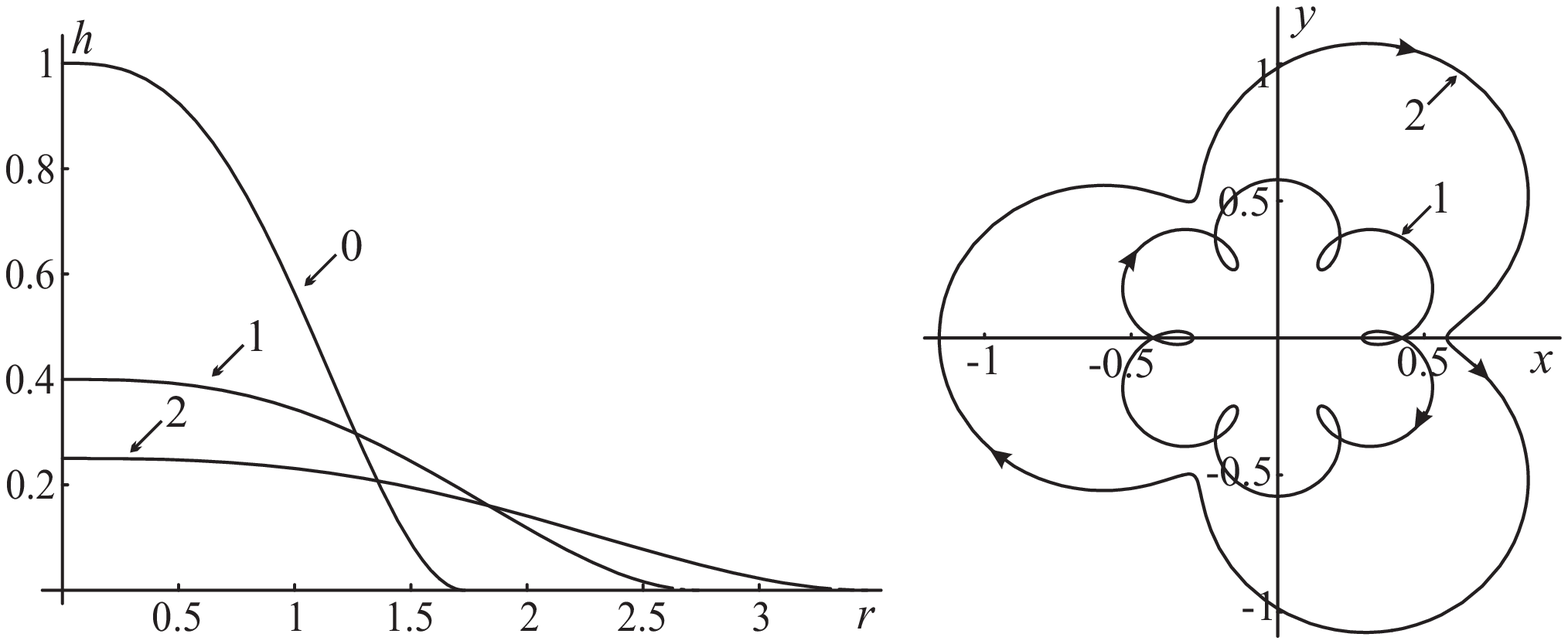}} \\[0pt]
\parbox{0.49\textwidth}{\caption{The fluid depth $h(t,r)$ at $t=\pi n/(2f)$, $n=0,1,2$
(lines 0, 1, and 2, correspondingly).} \label{fig:1}} \hfill
\parbox{0.49\textwidth}{\caption{Typical trajectories
of the fluid particles motion.} \label{fig:2}}
\end{center}
\end{figure}

Figure \ref{fig:1} presents the fluid depth depending on the radius $r$ at
different time moments $t=\pi n/(2f)$ with $n=0,1,2$ (lines 0, 1, and 2,
correspondingly). The fluid depth takes zero value at the circumference of
radius
\[ r=R_*(t)=
   -\frac{f}{l}\sqrt{\frac{1+\alpha^2\tau^2(t)}{\alpha(1+\tau^2(t))}}. \]
The graphs shown in figure \ref{fig:1} and others are plotted for $\alpha=2$,
$f=1$, $g=1$. As the parameters change, the qualitative view of the figure
remains the same.

Formulae (\ref{eq:tr-periodic}) imply that all particles, which belong to the
the circumference of radius $r_0=-\frac{f}{l\sqrt{\alpha}}\,\frac{m}{M}$ (here
$M \geq m>0$ are integers) at initial time $t=0$, have closed trajectories. At
$t=2\pi M/f$ the particles return to their starting position at $t=0$.
Generally, trajectories of fluid particles are quasi-closed. This means that
for any specified value $\varepsilon>0$ it is possible to choose
$t_{\varepsilon}>0$, so that the particle at $t=t_{\varepsilon}$ will be
located at a distance of no more than $\varepsilon$ from its initial position
at $t=0$. It should be pointed out that within the period of the solution
$T=2\pi/f$ only those particles, which belong to the boundary $r=R_*(t)$, where
fluid depth $h$ vanishes, return to their starting positions. In figure
\ref{fig:2} closed trajectories of particles, located at $t=0$ in the
circumference of radius $r=1/(2\sqrt{3})$ (curve 1) and $r=1/\sqrt{3}$ (curve
2), are shown. The graphics correspond to the choice $(m,M)=(1,6)$ and $(1,3)$,
with the above-indicated values $\alpha$, $f$ and $g$. Arrows show the
direction of the particles motion.

\begin{figure}[tbp]
\begin{center}
\resizebox{0.85\textwidth}{!}{\includegraphics{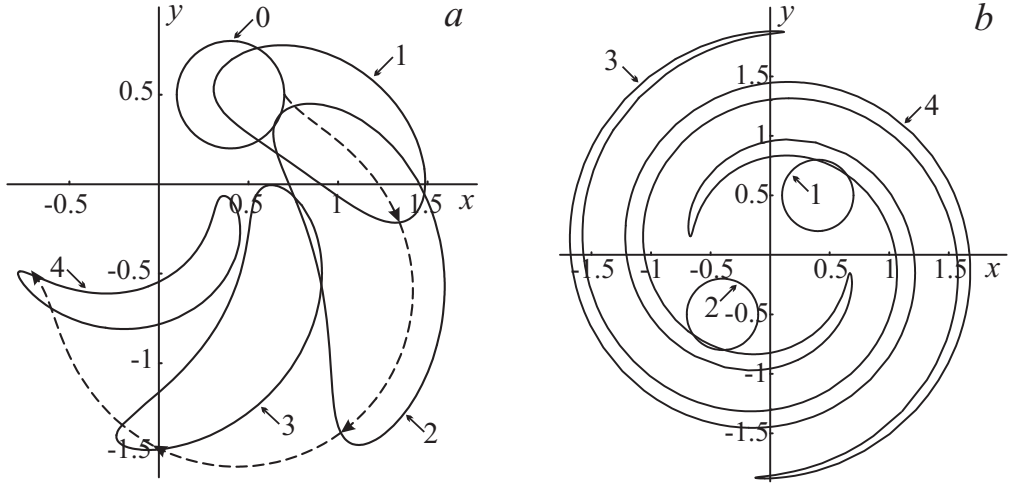}}\\[0pt]
\caption{The evolution of the material curves: $a$ --- short time behaviour;
$b$ --- long time behaviour.} \label{fig:3}
\end{center}
\end{figure}

Figure \ref{fig:3} shows the evolution of the material curves, which consist
always of the same fluid particles. In contrast to the previous solution
(\ref{eq:period-sol-1}), here an arbitrary closed material curve, except for
circumferences with centre at origin of coordinates, never returns to its
initial position. Propagation of the closed material curve (which is initially
the circumference of radius $r=0.3$ with its centre at point $(x,y)=(0.4,
0.5)$) is shown in figure \ref{fig:3}($a$) at the instants $t=\pi n/(2f)$ with
$n=0,...,4$ (curves 0--4). The trajectory of the fluid particle, initially
located at point $(x,y)=(0.7,0.5)$, is shown in figure \ref{fig:3}($a$) by
dashed line. As time growth this material curve transforms into helix. Figure
\ref{fig:3}($b$) presents two material curves which form circumferences at
$t=0$ (labeled as 1 and 2); and these material curves at the instant $t=7\pi/f$
(helixes 3 and 4, correspondingly).

\begin{remark}
The deformation of the material curve, shown in the Fig. \ref{fig:3}, is
typical for all solutions from the class (\ref{eq:periodic-sol}) where function
$\bar{V}$ is not identically zero. Indeed, we can integrate equations of
trajectories (\ref{eq:tr-eq}) and represent the solution in the form
(\ref{eq:tr-periodic}) with
$C=(r_0\sqrt{\alpha})^{-1}\bar{V}(r_0\sqrt{\alpha})$.
\end{remark}

\section{Rotationally-symmetric flows}

The infinitesimal symmetries (\ref{eq:oper-f-1})--(\ref{eq:oper-f-3}) allow us
to construct invariant and partially invariant solutions to the RSW model. This
section is devoted to a set of rotationally-symmetric submodels, which can be
reduced to ODEs. These submodels are derived by invariant reduction of the RSW
equations (\ref{eq:model-pol}) with the use of two-dimensional parameterized
classes of the optimal system of the subalgebra $\Theta L_9$. Therefore, the
obtained solutions are essentially different, that is one solution does not
reduce to the other by change of the variables. Note that these solutions can
be obtained from the corresponding invariant solutions of the SW equations
(\ref{eq:SW}) using transformation (\ref{eq:equiv-tr}).

\subsection{Stationary rotationally-symmetric flows in a ring}

Submodel
\begin{equation}
 \eqalign{
  UU_r-r^{-1}V^2-fV+gh_r=0, \cr
  U(V_r+r^{-1}V+f)=0, \quad (rUh)_r=0. }
 \label{eq:inv-mod-1}
\end{equation}
characterizes the class of stationary rotationally-symmetric solutions.
Equations (\ref{eq:inv-mod-1}) are derived according to the invariant reduction
of the RSW model (\ref{eq:model-pol}) with the use of subalgebra
$(Y_5,Y_7+Y_8)$ of admitted infinitesimal symmetries (in the polar coordinates
the generators have the following form $\hat{Y}_5=\partial_\theta$ and
$\hat{Y}_7+\hat{Y}_8=2f^{-1}\partial_t-\partial_\theta$).

Assuming $U=0$, integration of (\ref{eq:inv-mod-1}) yields (\ref{eq:st-sol-1}).
In case $U\neq 0$, we obtain the following result
\begin{equation}
   U=\frac{C_3}{rh}, \quad V=\frac{C_2}{r}-\frac{fr}{2}, \quad
   F(h,r)=h^3+\varphi_1(r)h^2+\varphi_2(r)=0,
  \label{eq:st-sol}
\end{equation}
where $C_i$ ($i=1,2,3$) are constants and
\[ \varphi_1(r)=\frac{1}{g}
   \left(\frac{f^2r^2}{8}+\frac{C_2^2}{2r^2}-C_1\right),
   \quad \varphi_2(r)=\frac{C_3^2}{2gr^2}. \]
The necessary condition of existence of solutions to the equation $F(h,r)=0$
with $h>0$ is inequality $C_1>|C_2|/2$, which guarantees that $\varphi_1<0$ for
some values $r$.

\begin{figure}[tbp]
\begin{center}
\resizebox{0.98\textwidth}{!}{\includegraphics{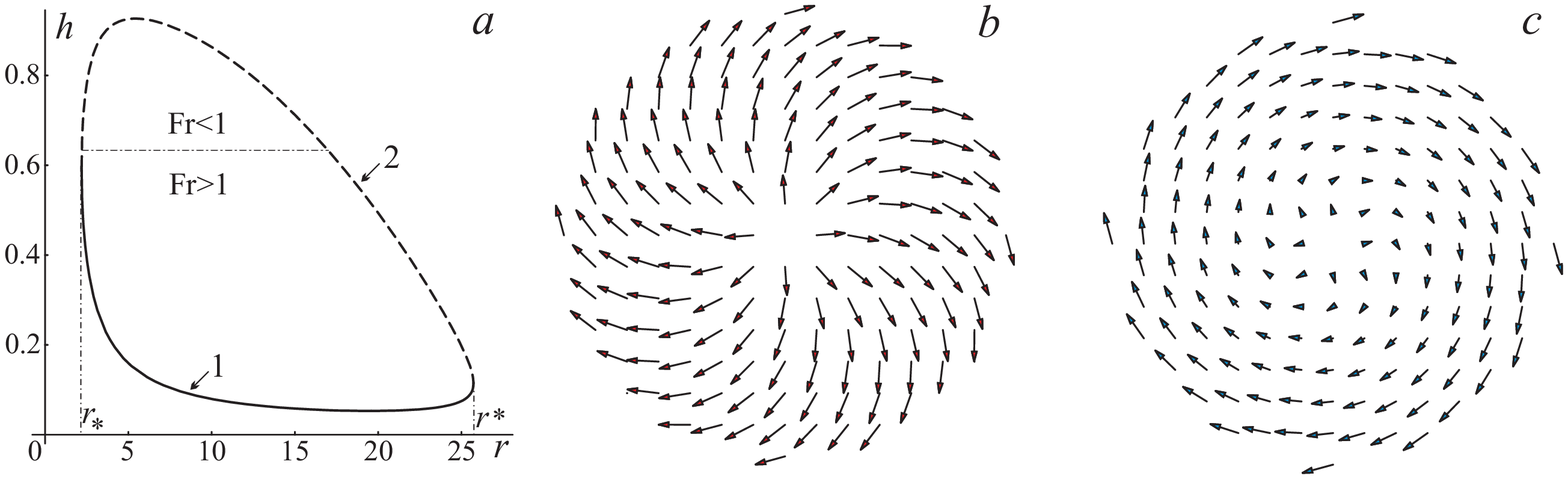}}\\[0pt]
\caption{Stationary rotationally-symmetric flow in a ring bounded by
characteristics: $a$ --- free surface height (line 1 corresponds to
supercritical flow, line 2 ---subcritical flow); $b$ --- velocity field
(supercritical flow); $c$ --- velocity field (subcritical flow).} \label{fig:4}
\end{center}
\end{figure}

The straightforward analysis shows that equation $F(h,r)=0$ has two branches of
solutions $h=h(r)$ defined within the interval $r_*\leq r \leq r^*$ (see figure
\ref{fig:4}($a$)). At the endpoints $r=r_*$ and $r=r^*$ the derivative $F_h$
vanishes: $F_h=(3h+2\varphi_1)h=0$. Hence, function $h(r)$ at the endpoints is
finite $h=h_c$ (here $h_c=-\frac{2}{3}\varphi_1(r_*)$ or
$h_c=-\frac{2}{3}\varphi_1(r^*)$), whereas its derivative is unbounded:
$h'(r_*)=h'(r^*)=-F_r/F_h=\infty$. The relation $U^2/2+gh=-g\varphi_1(r)$ is
fulfilled for the considered solution (\ref{eq:st-sol}), so that we get
equality $U^2=gh$ in circumferences of radii $r=r_*$ and $r=r^*$. Thus, the
surfaces $r=r_*$ and $r=r^*$ are sonic characteristics of the model. Recall
that the surface $\Phi(t,r,\theta)={\rm const}$ is the sonic characteristic of
equations (\ref{eq:model-pol}) on the solution $(U, V, h)$ if function
$\Phi(t,r,\theta)$ satisfies to equation
$\Phi_t+U\Phi_r+r^{-1}V\Phi_\theta=\pm\sqrt{gh}$.

The flow is supercritical (subcritical) in case ${\rm Fr}>1$ (${\rm Fr}<1$),
correspondingly. Here ${\rm Fr}=q/\sqrt{gh}$ is Froude number and
$q=\sqrt{U^2+V^2}$ is the magnitude of the velocity vector. It follows from
(\ref{eq:st-sol}) that $q^2=2C_1-C_2 f-2gh$. Hence, ${\rm Fr}=1$ at
$h=h_s=(2C_1-C_2 f)/(3g)$. One can easily prove that $h_c\leq h_s$. Therefore,
the lower branch of solution to equation $F(h,r)=0$ corresponds to the
supercritical flow. In figure \ref{fig:4}($a$) are presented the lower (line 1)
and the upper (line 2) branches of solution to equation $F(h,r)=0$. In figures
\ref{fig:4}($b$) and \ref{fig:4}($c$) the velocity fields $(u,v)$ are shown for
lower branch of free surface height $h=h(r)$ and upper branch, correspondingly.
The graphics are derived for $C_i=1$, $g=1$, $f=0.1$. As the parameters change,
the qualitative view of the graphics remains the same.

It is interesting to note that the class of stationary rotationally-symmetric
solutions of the RSW equations corresponds to the class of solutions of the SW
model which is invariant under rotation (generator $Z_5$) and nontrivial
projective transformation (generator $Z_7+Z_8$). On the basis of projective
transformation exact invariant solutions of the SW model and two-dimension gas
dynamic equations (with special state equation) were obtain and studied in
\cite{Nik, Gol04, Khab} and others.

\subsection{Collapse of a liquid ring. Regime \#1}

Let us consider the submodel derived by the invariant reduction of the RSW
equations (\ref{eq:model-pol}) using two-dimensional subalgebra $(Y_5,Y_7)$ of
infinitesimal transformations. According to the algorithm of invariant
solutions construction, let us solve the equations
\begin{equation*}
\hat{Y}_5 J=0, \quad \hat{Y}_7 J=0
\end{equation*}
to obtain the set of basic invariants
\begin{equation*}
  J=\left\{\frac{1-\cos(ft)}{r^2},\
  rU-\frac{fr^2}{2}\frac{\sin(ft)}{1-\cos(ft)},\ rV+\frac{fr^2}{2},\  r^2h
  \right\}.
\end{equation*}
The representation of an invariant solution is
\begin{equation}
  U=\frac{\varphi(\lambda)}{r}+\frac{fr}{2}\cot\Bigl(\frac{ft}{2}\Bigr), \quad
  V=\frac{\psi(\lambda)}{r}-\frac{fr}{2}, \quad h=\frac{\eta(\lambda)}{r^2},
 \label{eq:coll-1}
\end{equation}
where $\lambda=(1-\cos(ft))r^{-2}$. By substituting the representation of
solution into (\ref{eq:model-pol}), one obtains the system of ODEs for unknowns
$\varphi$, $\psi$, $\eta$
\begin{equation}
 (\lambda(\varphi^2+2g\eta))'+\psi^2=0, \quad \lambda\varphi\psi'=0, \quad
 (\lambda\varphi\eta)'=0.
 \label{eq:inv-mod-2}
\end{equation}

By virtue of the second equation (\ref{eq:inv-mod-2}) either $\varphi=0$ or
$\psi={\rm const}$. In the first case we get the following class of solutions
\begin{equation}
  \varphi=0, \quad \psi=\psi(\lambda), \quad \eta=\frac{1}{\lambda}
  \left(\lambda_0\eta_0-\frac{1}{2g}\int_{\lambda_0}^\lambda \psi^2(\nu)\,d\nu
  \right).
 \label{eq:coll-2}
\end{equation}
Here $\psi(\lambda)$ is an arbitrary smooth function, and $\eta_0$, $\lambda_0$
are constants. Solution (\ref{eq:coll-1}), (\ref{eq:coll-2}) is defined over
the finite time interval $t\in [t_0,2\pi/f)$, $t_0>0$. It has a singularity at
$r\to \infty$. In some cases, depending on the choice of function $\psi$, the
solution has singularity at $r\to 0$.

Notice, that surfaces $\lambda(t,r)={\rm const}$ are contact characteristics on
the solution (\ref{eq:coll-1}), (\ref{eq:coll-2}). This allows one to interpret
the solution as a collapse of a liquid ring (or a liquid cylinder) $R_1(t)\leq
r \leq R_2(t)$ compressed by pistons which move according to the law
\begin{equation*}
R_i(t)=R_{i0}\sin(ft/2) \quad (0\leq R_{10}<R_{20}),
\end{equation*}
so that the pistons surfaces coincide with the contact characteristics. In the
course of the collapse of a liquid ring (with $t\to 2\pi/f$), the radial
velocity $U=\frac{fr}{2}\cot\left(\frac{ft}{2}\right)$ in the range of
$R_1(t)\leq r \leq R_2(t)$ is bounded, whereas the liquid depth $h$ as well as
the circular velocity $V$ (with any $\psi\neq {\rm const}$) increases
infinitely.

The other class of solutions of submodel (\ref{eq:inv-mod-1}) is derived in
case $\varphi \neq 0$, $\psi=C_2={\rm const}$. Whereas,
$\eta(\lambda)=C_3/(\lambda\varphi)$ and function $\varphi(\lambda)$ are
defined from cubic equation
$\varphi^3+(C_2^2-C_1/\lambda)\varphi+2gC_2/\lambda=0$. The solution of the RSW
equations (\ref{eq:model-pol}) is obtained by using representation of solution
(\ref{eq:coll-1}).

These unsteady solutions of the RSW equations can be also derived from
well-known stationary rotational-symmetric solutions of the SW equations using
transformations (\ref{eq:equiv-tr}) which relate the models.

\subsection{Collapse of a liquid ring. Regime \#2}

Submodel
\begin{equation}
   \varphi'=(\psi+f)\psi-\varphi^2-2g\eta, \quad \psi'=-(2\psi+f)\varphi, \quad
   \eta'=-4\varphi\eta.
 \label{eq:inv-mod-3}
\end{equation}
arises from the invariant reduction of the RSW equations (\ref{eq:model-pol})
with the use of the subalgebra $(Y_5,Y_6)$. Indeed, in polar coordinates
(\ref{eq:polar}) generator $Y_6$ has the following form:
$\hat{Y}_6=r\partial_r+U\partial_U+V\partial_V+2h\partial_h$. The set of basic
invariants of the generators $\hat{Y}_5$ and $\hat{Y}_6$ is
\begin{equation*}
J=\{t, r^{-1}U,r^{-1}V, r^{-2}h\}.
\end{equation*}
Thus, we get the following representation of the invariant solution
\begin{equation}
  U=r\varphi(t), \quad V=r\psi(t), \quad h=r^2\eta(t).
 \label{eq:ring-1}
\end{equation}

Substitution of representation (\ref{eq:ring-1}) into the RSW model
(\ref{eq:model-pol}) produces ODEs (\ref{eq:inv-mod-3}). In the general case it
is difficult to get a solution of the submodel (\ref{eq:inv-mod-3}) in closed
form. Let us take $\psi(t)=-f/2$. Then the second equation in
(\ref{eq:inv-mod-3}) is fulfilled automatically, whereas the solution to the
first and the third equations is given in an implicit form
\begin{equation}
 \fl \varphi=\hat{\varphi}(\eta)=
  \pm\sqrt{2g\eta-\frac{f^2}{4}+\left(\varphi_0^2-2g\eta_0+
  \frac{f^2}{4}\right)\sqrt{\frac{\eta}{\eta_0}}}, \quad
  t(\eta)=-\frac{1}{4}\int_{\eta_0}^\eta \frac{d\,\nu}{\nu\hat{\varphi}(\nu)}.
  \label{eq:ring-2}
\end{equation}
Here arbitrary constants $\varphi_0$ and $\eta_0>0$ are values of functions
$\varphi(t)$ and $\eta(t)$ at $t=0$ (initial data for ODEs
(\ref{eq:inv-mod-3})). Function $\varphi(t)$ is strictly decreasing by virtue
of the first equation (\ref{eq:inv-mod-3}) (due to $\psi=-f/2$). It follows
from the last equation (\ref{eq:inv-mod-3}) that function $\eta(t)$ increases
as time grows if $\varphi(t)<0$, otherwise $\eta(t)$ decreases.

If constant $\varphi_0\leq 0$, the minus sign precedes the root in the first
formula (\ref{eq:ring-2}). In fact, $\varphi(0)=\hat{\varphi}(\eta_0)\leq 0$
and $\varphi(t)$ is strictly decreasing function. It is evident from the second
formula (\ref{eq:ring-2}) that $t\to T_*>0$ as $\eta\to \infty$ (integral in
(\ref{eq:ring-2}) converges). The solution (\ref{eq:ring-1}), (\ref{eq:ring-2})
can be interpreted as collapse of a liquid ring during the finite time interval
$[0,T_*]$. At $t=0$ liquid is located in the ring (or in the cylinder) $0\leq
R_{10}\leq r\leq R_{20}$. As time grows, impermeable boundaries of the ring
move in accordance with the law
\[ r= R_i(t)=\hat{R}_i(\eta)=R_{i0}(\eta_0/\eta)^{1/4}, \]
providing the fulfillment of boundary condition $R_i'(t)=U(t,R_i(t))$. In
contrast to the previous solution (\ref{eq:coll-1}), (\ref{eq:coll-2}) at $t\to
T_*$ the circular velocity $V$ vanishes, whereas the radial velocity $U$ and
the depth $h$ infinitely increase. At that, the moving pistons collapse into
the origin $r=0$.

Let us take $\varphi_0>0$. In view of
$\varphi_0=\varphi(0)=\hat{\varphi}(\eta_0)>0$, the plus sign precedes the root
in first formula (\ref{eq:ring-2}). As mentioned above, function $\eta$
decreases while $\varphi>0$. It is obvious from (\ref{eq:ring-2}) that radicand
vanishes at some value $\eta=\eta_1<\eta_0$. At time $t_1=t(\eta_1)<\infty$
function $\varphi(t)$ changes sign from positive to negative and function
$\eta(t)$ grows for $t>t_1$. The behaviour of solution for $t>t_1$ is similar
to the previous case $\varphi_0\leq 0$. Thus, the solution can be interpreted
as a partial spreading $(0\leq t\leq t_1)$ and collapse $(t_1<t<T_*)$ of a
liquid ring.

\section{Conclusions}

In this paper we have studied symmetry properties and some classes of exact
solutions of the RSW model using group analysis. We have determined the
9-dimensional Lie algebra of admissible infinitesimal generators. We have shown
that the derived Lie algebra of symmetries is isomorphic to the Lie algebra of
infinitesimal transformations admitted by (2+1)-dimensional SW equations. This
allows one to employ its known optimal system of subalgebras for construction
of essentially different invariant solutions of the RSW equations. Moreover, we
have found change of variables which reduce the RSW equations to the SW model.
This transformation allow one to construct and study solutions of the RSW
equations using solutions of the SW model and vice versa. We have also derived
and analysed finite transformations corresponding to the nontrivial symmetries
of the RSW model. Using these transformations we have generated new
time-periodic exact solutions of the RSW equations. These solutions describe
fluid flow with quasi-closed particle trajectories and may be interpreted as
pulsation of liquid volume under the influence of gravity and Coriolis forces.
Based on some two-dimensional parameterized classes of the optimal system of
subalgebras we have reduced the RSW model to the ODEs and integrate them. In
particular, we have constructed and studied exact solutions of the RSW model
describing rotational-symmetric stationary flows in a ring bounded by
characteristics as well as time dependent flows which describe spreading and
collapse of a liquid ring. We have also pointed out exact solutions of the SW
equations corresponding to the obtained solutions of the RSW equations.

\ack The work was supported by the Russian Foundation for Basic Research
(project 07-01-00609), Programme of support of the Leading scientific schools
(grant 2826.2008.1) and by the Russian Academy of Sciences (project 4.14.1).

\end{document}